\documentclass[final]{siamart171218}

\usepackage[english]{babel}
\usepackage[utf8x]{inputenc}
\usepackage[T1]{fontenc}
\usepackage{theorem}
\usepackage{secdot}
\usepackage{enumerate}
 
\usepackage{amsmath,amssymb}
\usepackage{graphicx}
\usepackage{subfig}
\usepackage[colorinlistoftodos]{todonotes}

\newcommand{\expect}{\mathbb{E}}

\DeclareMathOperator*{\argmin}{\arg\min}

\title{Preconditioning Kaczmarz method by sketching\thanks{This research was supported by the RFBR grant 18-31-20069-mol\_a\_ved}}
\author{Alexandr Katrutsa\footnotemark[2]\ \footnotemark[3] \and Ivan Oseledets\footnotemark[2]}

\headers{Preconditioning Kaczmarz method by sketching}{Alexandr Katrutsa and Ivan Oseledets}

\begin{document}
\maketitle

\renewcommand{\thefootnote}{\fnsymbol{footnote}}
\footnotetext[2]{Skolkovo Institute of Science and Technology, Nobel St., 3, 143025, Moscow, Russia, (\email{aleksandr.katrutsa@phystech.edu}, \email{i.oseledets@skoltech.ru})}
\footnotetext[3]{Moscow Institute of Physics and Technology, Institutskii per. 9, Dolgoprudny, 141700, Moscow Region, Russia}
\begin{abstract}
We propose a new method for preconditioning Kaczmarz method by sketching.
Kaczmarz method is a stochastic method for solving overdetermined linear systems based on a sampling of matrix rows. 
The standard approach to speed up convergence of iterative methods is using preconditioner.
As we show the best possible preconditioner for this method can be constructed from QR decomposition of the system matrix, but the complexity of this procedure is too high.
Therefore, to reduce this complexity, we use random sketching and compare it with the Kaczmarz method without preconditioning.
The developed method is applicable for different modifications of classical Kaczmarz method that were proposed recently.
We provide numerical experiments to show the performance of the developed methods on solving both random and real overdetermined linear systems.  
% This paper studies different approaches to speed up convergence of Kaczmarz method using preconditioners.
% Preconditioning is well-known approach based on reduction condition number of the matrix in linear system. 
% We employ this approach to speed up Kaczmarz method since existing convergence results claim dependence of convergence speed on condition number of the matrix.  
% Kaczmarz method is the method to solve overdetermined linear systems.
% These systems also can be solved in the least-squares manner, therefore we compare Kaczmarz method with other stochastic methods that solve least squares problem associated with considered linear system.
\end{abstract}



\section{Introduction}
This paper presents a new method of preconditioning of the Kaczmarz method~\cite{kaczmarz1379} based on sketching technique~\cite{woodruff2014sketching}.
The Kaczmarz method is the stochastic method for solving overdetermined linear systems $Ax = b$, where $A \in \mathbb{R}^{m \times n}$ and $m \gg n$.
We assume that we can only sample rows of $A$, and can not use the~full~$A$. 
Existing convergence results~\cite{strohmer2009randomized} show that convergence is defined by $ \kappa_F(A) = \|A\|_F\|A^{-1}\|_2$.
We show by extensive numerical simulations that this value can be minimized by using preconditioner.
To reduce the time required for a preconditioner computation we approximate the matrix~$A$ by sketching and compute the preconditioner for a sketched matrix.
According to our assumptions, we can sketch the matrix only by random sampling of a fixed number of rows.
The number of sketched rows is a parameter of the proposed method, and we test different numbers proportional to the number of columns $n$.
% Unfortunately, study~\cite{avron2010blendenpik} shows that this method of sketching can give very bad preconditioner and proposes methods to fix this issue.
% However, these methods are too costly and can not be applicable in our case.

Numerical experiments demonstrate the performance of the proposed method for both random consistent linear systems, test tomography problems and non-linear function approximation problem with random Fourier features technique.
Also, we study random linear systems with noisy right-hand sides and show how the variance of the noise affects the convergence.

Main contributions of this paper are the following.
\begin{itemize}
    \item We study the preconditioning of the Kaczmarz method by sketching, analyze its complexity and provide limitations of this approach. (Section~\ref{sec::kacz},~\ref{sec::sketch}).
    \item We propose to use preconditioning as the fine-tuning technique to get highly accurate solutions of overdetermined linear systems. (Section~\ref{sec::num_exp})
    \item We compare the proposed approach on both synthetic data, test tomography data and function approximation problem using explicit feature map. (Section~\ref{sec::num_exp})
\end{itemize}

\subsection{Related work}
The Kaczmarz method for solving overdetermined linear systems is a well-known randomized iterative method, which has numerous interpretations and generalizations~\cite{gower2015randomized,needell2014paved}.
For example, it can be re-written as a special case of the stochastic gradient descent or coordinate descent applied to the linear least-squares problem corresponding to the considered system~\cite{needell2014stochastic,qu2016coordinate}.
The idea of preconditioning of such kind of methods is mentioned in the paper~\cite{yang2017weighted} but in the context of using stochastic gradient method for solving $\ell_p$ regression problem.
In contrast, we emphasize on the Kaczmarz method and its convergence speed up with preconditioning.
The convergence of the Kaczmarz method for a noisy linear system is presented in~\cite{needell2010randomized}.
Our study shows that in the case of small noise the preconditioned Kaczmarz method gives smaller relative error than without preconditioning.
Also, there is an extension~\cite{zouzias2013randomized} of the Kaczmarz method for the inconsistent linear system that uses random orthogonal projection technique to address the inconsistency.
One of the main issues in the Kaczmarz method is the correct sampling of rows: if two sequential rows are correlated than the progress in the convergence will be very small. 
The paper~\cite{needell2013two} addresses this issue and proposes a two-subspace projection method.

\section{Kaczmarz method and its preconditioning}
\label{sec::kacz}
Let $A \in \mathbb{R}^{m \times n}$ be a rectangular matrix, where $m \gg n$ and $b \in \mathbb{R}^m$. 
We consider the problem of solving an overdetermined linear system
\begin{equation}
    Ax = b,
    \label{eq::ovdet_sys}
\end{equation}
with the following assumptions
\begin{enumerate}
    \item[1)] the least-squares solution of the problem using a full matrix $A$ is not available
    \item[2)] we can sample rows of matrix $A$ according to some probability distribution. 
\end{enumerate}
To fulfill these assumptions, we consider the Kaczmarz method~\cite{kaczmarz1379} for solving~\eqref{eq::ovdet_sys}. 
This method is based on the iterative update of the current approximation~$x_k$ using the solution of the following optimization problem:
\begin{equation}
\begin{split}
& x_{k+1} = \argmin_x \frac{1}{2}\|x_k - x\|^2_2, \\
\text{s.t. } & a_ix = b_i,
\end{split}
\label{eq::kaczmarz_problem}
\end{equation}
where $a_i$ is the $i$-th row of the matrix $A$ and $b_i$ is the $i$-th element of the vector $b$.
Problem~\eqref{eq::kaczmarz_problem} has an analytical solution
\begin{equation}
x_{k+1} = x_k - \frac{a_ix_k - b_i}{\|a_i\|_2^2}a^{\top}_i.
\label{eq::kaczmarz_update}
\end{equation}
The rows of the matrix $A$ can be sampled in different ways. 
The original paper~\cite{kaczmarz1379} proposed a sequential selection of rows one by one.
Later, the study~\cite{strohmer2009randomized} proposed a \emph{Randomized Kaczmarz method} (RKM), which samples rows with probabilities $\frac{\|a_i\|_2^2}{\|A\|_F^2}$.
The convergence of the RKM is characterized by the following Theorem.

 \begin{theorem}[\cite{strohmer2009randomized}]
Let $x^*$ be the solution of~\eqref{eq::ovdet_sys}.
Then RKM converges in expectation with the average error as
\[
\expect \|x^* - x_k \|_2^2 \leq (1 - \kappa_F(A)^{-2})^k \| x^* - x_0\|^2_2,
\]
where $\kappa_F(A) = \|A\|_F \|A^{-1}\|_2 = \frac{\| A \|_F}{\sigma_{\min}(A)}$.
\label{th::rka_ub}
\end{theorem}

From this Theorem, it follows that decreasing $\kappa_F(A)$ leads to faster convergence of RKM.
Also, let $k(A) = \|A\|_2\|A^{-1}\|_2$ be a condition number of the matrix $A$. 
The following well-known inequality between $\kappa_F(A)$ and $k(A)$ holds
% \todo[inline]{Something about tightness of inequality?}
\begin{equation}
\kappa_F(A) \leq \sqrt{n} \cdot k(A).
\label{eq::k_vs_kappa}
\end{equation}
Thus, $k(A)$ is an upper bound of $\kappa_F(A)$ up to the factor $\sqrt{n}$, and therefore instead of the minimization of $\kappa_F(A)$ we can minimize~$k(A)$.
The standard approach for the reduction of the condition number is to use \emph{preconditioners}~\cite{saad2003iterative}.
We will denote by $P_L \in \mathbb{R}^{m \times m}$ and by $P_R \in \mathbb{R}^{n \times n}$ the left and right preconditioners, respectively.
In our case, the matrix $A$ is rectangular of size $m \times n$ and $m \gg n$.
Therefore, the right preconditioner is preferable since it requires $\mathcal{O}(n^2)$ memory in contrast to the left preconditioner, which requires~$\mathcal{O}(m^2)$. 
An ideal right preconditioner $P_R^*$ satisfies 
\[
k(AP_R^*) = 1,
\]
which means that $AP_R^*$ is unitary.
Therefore, we can obtain $P_R^*$ from the QR decomposition of the matrix $A$. 
In particular, if 
\begin{equation}
    A = QR,
    \label{eq::qr_full}
\end{equation}
where $Q \in \mathbb{R}^{m \times n}$ is an orthogonal matrix and $R \in \mathbb{R}^{n \times n}$ is an upper triangular matrix, then 
\begin{equation}
    Q = AR^{-1}, \quad P_R^* = R^{-1}.
    \label{eq::opt_prec}
\end{equation}
The complexity of every step in this procedure for deriving optimal right preconditioner $P^*_R$ is the following:
\begin{itemize}
    \item QR decomposition of the matrix $A \in \mathbb{R}^{m \times n}$ costs $\mathcal{O}(mn^2)$,
    \item inversion of upper triangular matrix $R$ costs $\mathcal{O}(n^3)$.
\end{itemize}
Thus, the total complexity is $\mathcal{O}(mn^2) + \mathcal{O}(n^3)$, which is too high in the case of $m \gg n$. 
Section~\ref{sec::sketch} shows how to reduce the complexity of preconditioner computation down to $\mathcal{O}(n^3)$.

After introducing the right preconditioner $P_R$ we have the following modification of the problem~\eqref{eq::ovdet_sys}:
\begin{equation}
    AP_Ry = b, \quad x = P_Ry.
\end{equation}
The corresponding modification of the Kaczmarz method is
\begin{equation}
    x_{k+1} = x_k - \frac{a_iP_Rx_k - b_i}{\|P_Ra_i\|_2^2}(a_iP_R)^{\top},
    \label{eq::prekacz_update}
\end{equation}
where $a_i$ is the $i$-th row of the matrix $A$ sampled according to some distribution.

Note that the cost of the one iteration of the method~\eqref{eq::prekacz_update} is $\mathcal{O}(m^2 + n)$, which is higher than corresponding cost $\mathcal{O}(n)$ of the original Kaczmarz method~\eqref{eq::kaczmarz_update}. 
Also, we need to perform preprocessing for computing $P_R$.
However, the expected speedup asymptotically dominates the increased cost of every iteration and gives a smaller relative error for the same time.
Examples of such behavior and runtime comparison are provided in Section~\ref{sec::num_exp}.

Although we show in~\eqref{eq::qr_full} and~\eqref{eq::opt_prec} how to compute the optimal right preconditioner, this computation is intractable because we can not use the full matrix $A$.
To overcome this limitation we propose to use \emph{sketching}~\cite{woodruff2014sketching,ghashami2016frequent,gower2018stochastic} of the matrix $A$ to compute the approximation of $P^*_R$.
The next section discusses some ways to sketch the matrix $A$ and resulting sketching-based preconditioners.

\section{Sketching-based preconditioners}
\label{sec::sketch}
This section gives an overview of \\ sketching idea, how it helps in constructing the preconditioner for the Kaczmarz method and limitations of this approach.

The sketching is the procedure of replacing original matrix $A \in \mathbb{R}^{m \times n}$ by a sketched matrix $\hat{A} = SA \in \mathbb{R}^{r \times n}$, where $S \in \mathbb{R}^{r \times m}$ is some given matrix and \\ $r \ll m$~\cite{woodruff2014sketching,gower2015randomized}.
Since the complexity of computing the product $SA$ for arbitrary dense matrices $S$ and $A$ is $\mathcal{O}(rmn)$, which is intractable, we consider matrix~$S$ from a set~$\mathcal{S}$, such that $S = [s_{ij}] \in \mathcal{S} \subset \mathbb{R}^{r \times n}$ if and only if for $i = 1,\ldots,r$, 
\[
s_{ij} = 
\begin{cases}
1 & j = k_i, \\
0 & \text{otherwise},
\end{cases}
\]
where $k_1 < k_2 < \ldots < k_r$ and $\{k_1, k_2, \ldots, k_r\} \subset \{1, \ldots, m\}$.
For such matrices $S \in \mathcal{S}$, product $SA$ is equivalent to the selection of $r$ rows of the matrix $A$ with indices $\{k_1, \ldots, k_r\}$ such that $k_1 < k_2 < \ldots < k_r$.
Details of rows selection we discuss later in this section.
Now assume we have a sketched matrix $\hat{A}$ such that $r \geq n$. 
Then we compute QR decomposition of this matrix 
\[
\hat{A} = \hat{Q}\hat{R},
\]
where $\hat{Q} \in \mathbb{R}^{r \times n}$ is orthogonal and $\hat{R} \in \mathbb{R}^{n \times n}$ is upper triangular matrix,
and use $\hat{P}_R = \hat{R}^{-1}$ as an approximation of $P^*_R$. 
Comparing with~\eqref{eq::qr_full} and~\eqref{eq::opt_prec}, we reduce the complexity of the QR decomposition from $\mathcal{O}(mn^2)$ to $\mathcal{O}(rn^2)$, which can be rewritten in the form $\mathcal{O}(n^3)$ if $r = \gamma n$, where $\gamma \geq 1$ is some given factor.
Therefore, the total complexity of computing preconditioner $\hat{P}_R$ is $\mathcal{O}(n^3)$ and does not depend on the number of rows $m$ of the matrix $A$.  

\subsection{Selection of rows for sketched matrix $\hat{A}$}

The remaining question in the procedure of construction~$\hat{P}_R$ is how to obtain a sketched matrix~$\hat{A}$ or in particular how to select rows in the matrix~$A$ in such a way that the matrix~$\hat{P}_R$ approximates the~matrix~$P^*_R$.

Since we can only sample rows of matrix $A$, but can not use full matrix, we use a simple uniform random sampling strategy. 
Unfortunately, this strategy is not universal and can give a very poor approximation of $P^*_R$ or even singular matrix $\hat{R}$.
This phenomenon relates to the notion of \emph{coherence} of the matrix.

\begin{definition}[\cite{avron2010blendenpik}]
Let $\mu(A)$ be the coherence of matrix $A$:
\begin{equation}
    \mu(A) = \max_{i=1,...,m} \|q_i \|^2_2,
\end{equation}
where $q_i$ is the $i$-th row of the matrix~$Q$ from the QR decomposition of the matrix~$A$.
\end{definition}
According to~\cite{avron2010blendenpik}, the coherence $\mu(A)$ is between $ n / m$ and $1$ and the larger coherence of the matrix, the poorer approximation $\hat{R}$ is given by a uniform random sampled sketched matrix~$\hat{A}$.
To address this issue, authors of~\cite{avron2010blendenpik} propose several row mixing methods that reduce the coherence of the matrix with unitary transformations, e.g, Discrete Cosine Transform.
Unfortunately, such kind of methods require $\mathcal{O}(m\log m)$ operations and are too costly, since in our case $m \gg n$.
Thus, the used approach of constructing sketched matrix~$\hat{A}$ is applicable only to matrices that have low coherence.

\section{Numerical experiments}
\label{sec::num_exp}
In this section, we provide numerical experiments to demonstrate the performance of the proposed approach to speed-up the convergence of the Kaczmarz method.
The source code can be found at GitHub\footnote{\url{https://github.com/amkatrutsa/preckacz}}.
We use random matrices, standard tomography data from AIRTool~II~\cite{hansen2018air} and matrices taken from non-linear function approximation problem.
We reduce the latter problem to the problem of solving the overdetermined linear system with random Fourier features technique~\cite{ambikasaran2016fast}.
In all experiments, we initialize $x_0$ with zero vector.
Considered algorithms are implemented in Python with NumPy library~\cite{oliphant2006guide}. 
Since we assume that the full matrix $A$ is not available, we use uniform random sampling of rows in the Kaczmarz method and its preconditioning modifications.
Also, we test the proposed approach for a range of factors $\gamma$ required to compute~$\hat{P}_R$.
We take into account the preprocessing time required to compute $\hat{P}_R$ as a starting point for all convergence plots shown below.

\subsection{Synthetic data}
\label{sec::random_data}
We generate a random matrix of size $m = 10000, n = 100$ with condition number $k \sim 10^6$.
Entries of the matrix are generated from the standard normal distribution $\mathcal{N}(0, 1)$.
After that, we generate a ground truth solution~$x^* \in \mathbb{R}^n$ from $\mathcal{N}(0, I)$.
We compare the proposed preconditioning method in two settings: consistent system, where $b = Ax^*$ and consistent system with noisy right-hand side: $b = Ax^* + \varepsilon$, $\varepsilon \sim \mathcal{N}(0, \sigma^2 I)$.
In both settings, we test different values of $\gamma$.
On the one hand, the larger $\gamma$, the better approximation quality of $\hat{P}_R$.
But on the other hand, the larger $\gamma$, the more time we need to compute $\hat{P}_R$ in pre-processing step.
Therefore, we identify the minimal value of $\gamma$ among the tested values that give faster convergence than the Kaczmarz method without preconditioning (``No preconditioning'' label in figures below).
% {\color{red}{Experiments with randomly generated overdetermined linear systems, which are consistent and consistent with noise, show which values of $\gamma$ are acceptable for faster convergence compared with the Kaczmarz method (``No preconditioning'' label in figures below), and in the same time do not lead to an overhead in the pre-processing step for computing $\hat{P}_R$. }}

\subsubsection{Consistent system setting}
\begin{figure}[!h]
    \centering
    \includegraphics[scale=0.3]{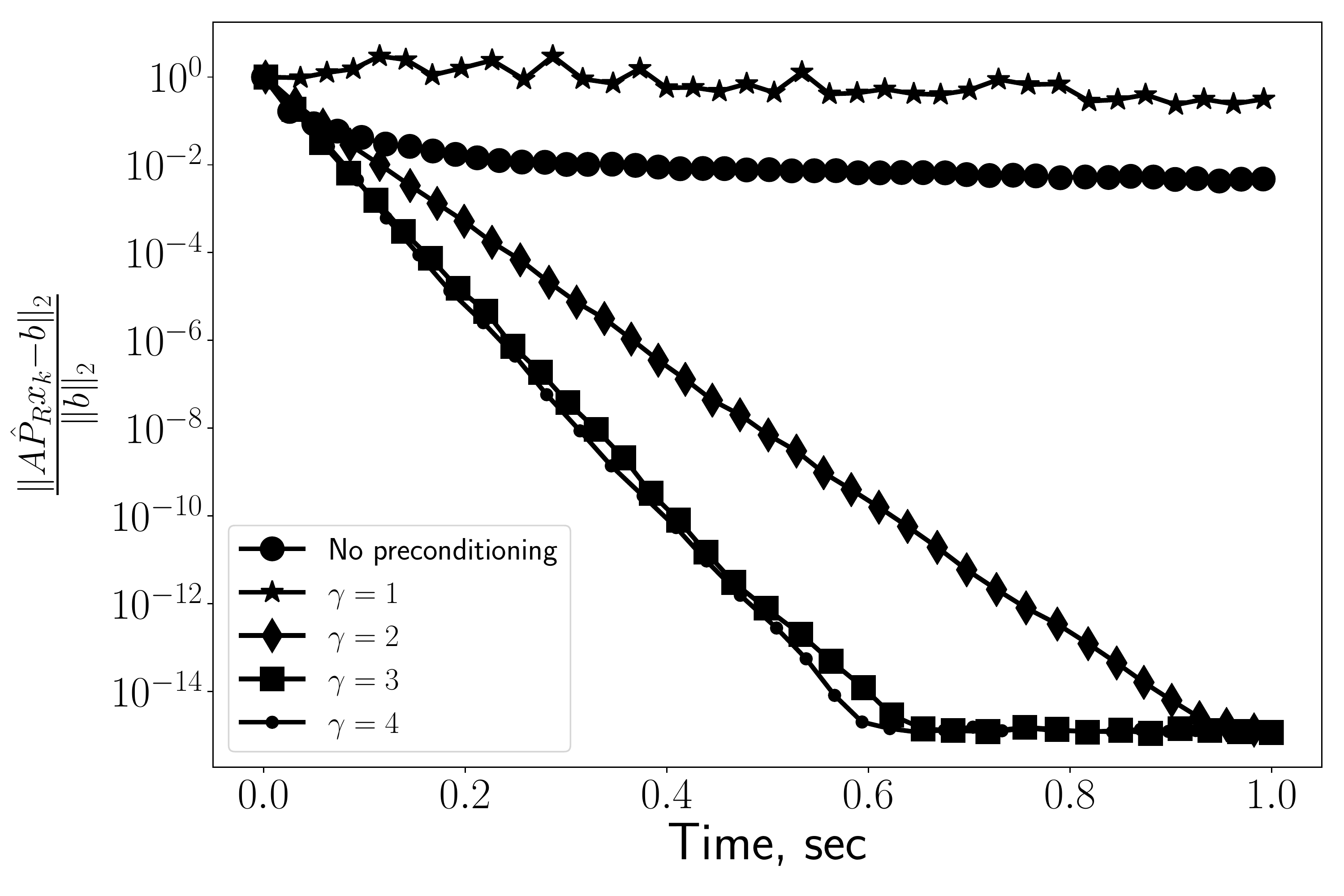}
    \caption{Comparison of convergence speed of the considered methods for consistent linear system.}
    \label{fig::random_consistent_time}
\end{figure}

In this experiment setting we generate the right-hand side $b = Ax^*$. 
Figure~\ref{fig::random_consistent_time} shows the dependence of the relative error on the convergence time for the compared methods. 
In the case of the Kaczmarz method, which is labeled as ``No preconditioning'', $\hat{P}_R = I$.
This plot demonstrates that if $\gamma$ is not large enough, e.g. $\gamma = 1$, the convergence of the preconditioned Kaczmarz method can be significantly slower than the convergence without preconditioning.
In the same time, already for $\gamma = 2$ we observe a significant speedup compared with the Kaczmarz method without preconditioning.
Further increase of $\gamma$ shows that $\gamma = 3$ is already large enough and a further increase is not reasonable.

\subsubsection{Consistent system with noise setting}

\begin{figure}[!h]
    \centering
    \subfloat[$\sigma = 10^{-1}$]{
    \centering
    \includegraphics[scale=0.2]{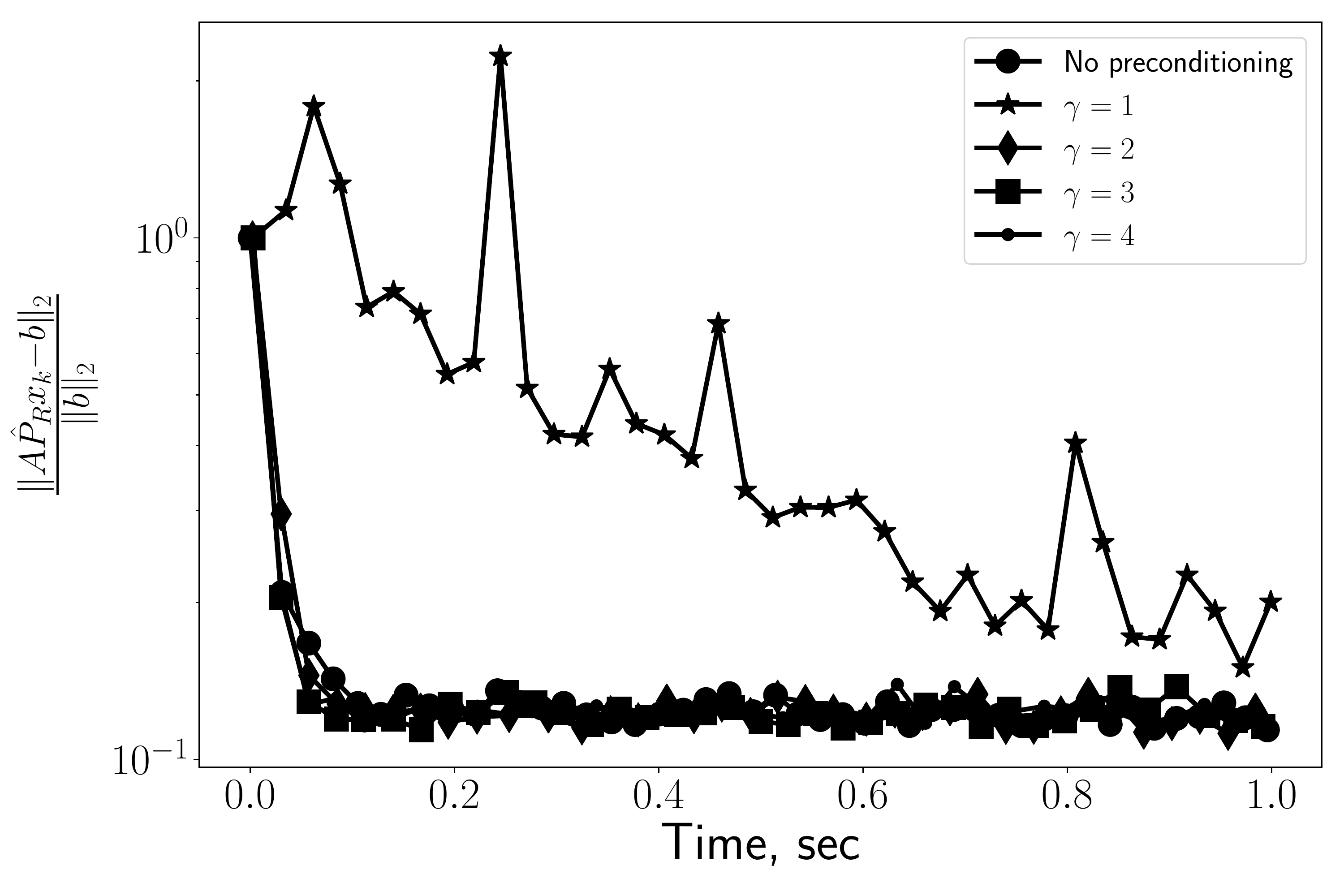}}
    \subfloat[$\sigma = 10^{-2}$]{
    \centering
    \includegraphics[scale=0.2]{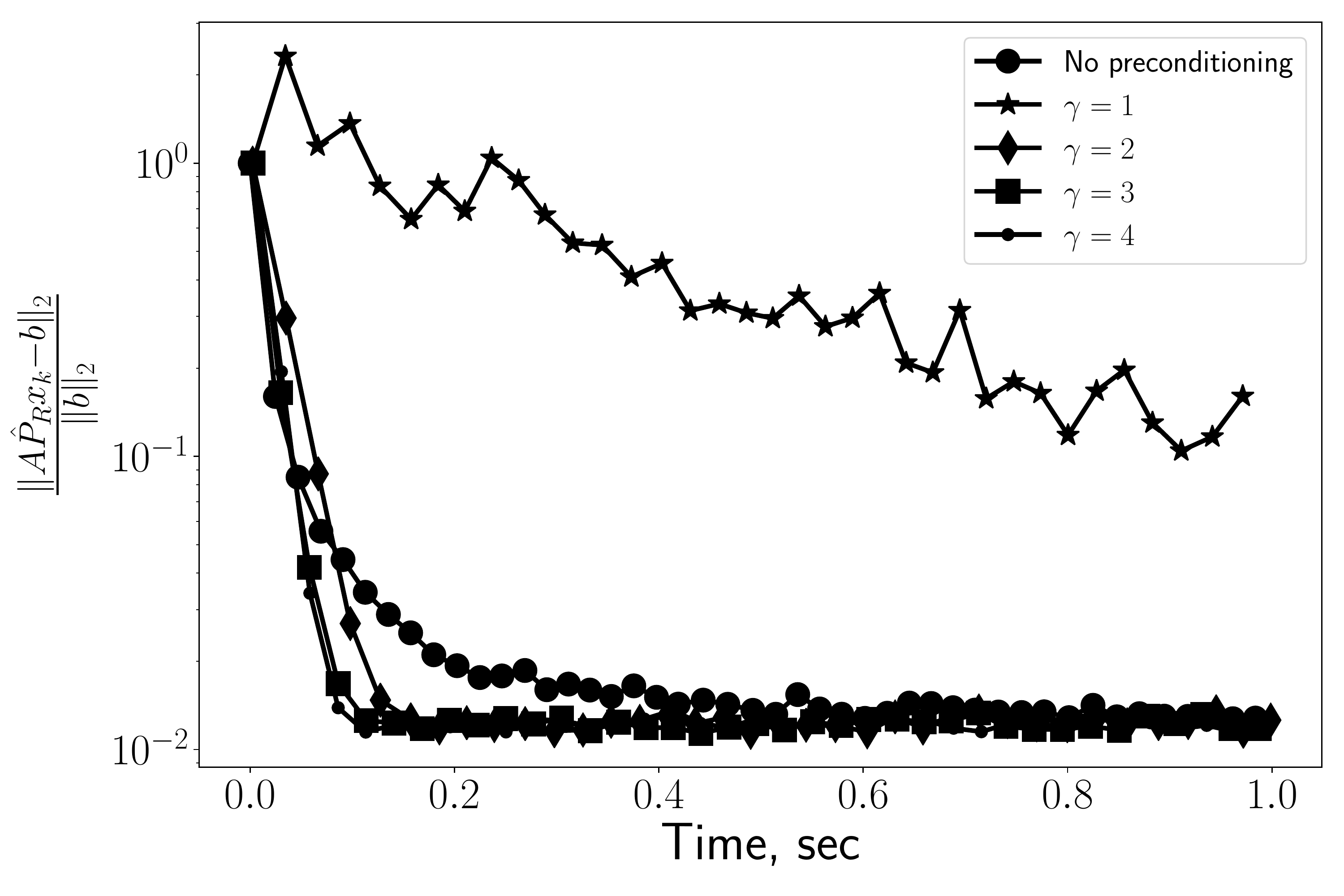}}
    \\
    \subfloat[$\sigma = 10^{-3}$]{\centering\includegraphics[scale=0.2]{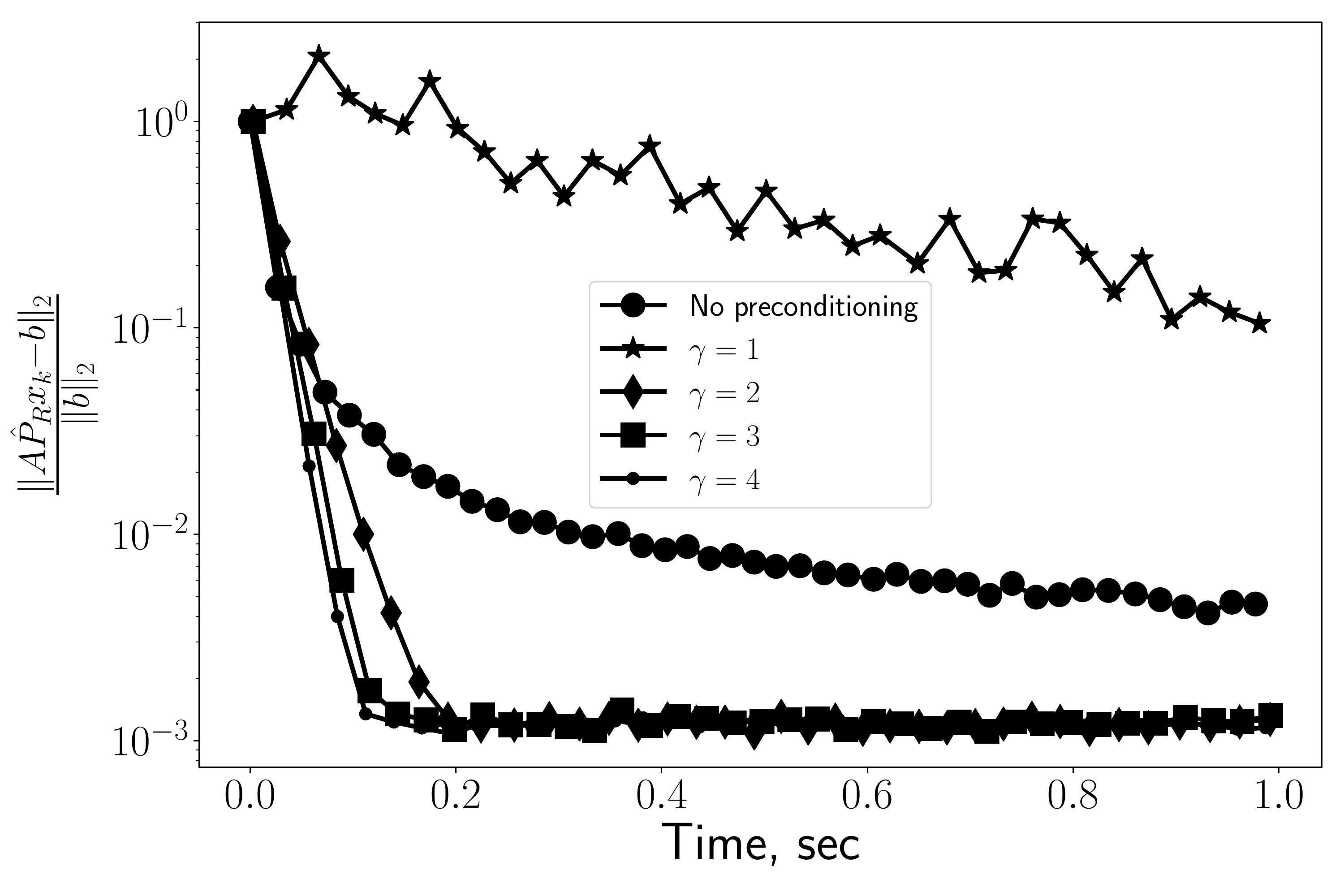}}
    \subfloat[$\sigma = 10^{-4}$]{\centering \includegraphics[scale=0.2]{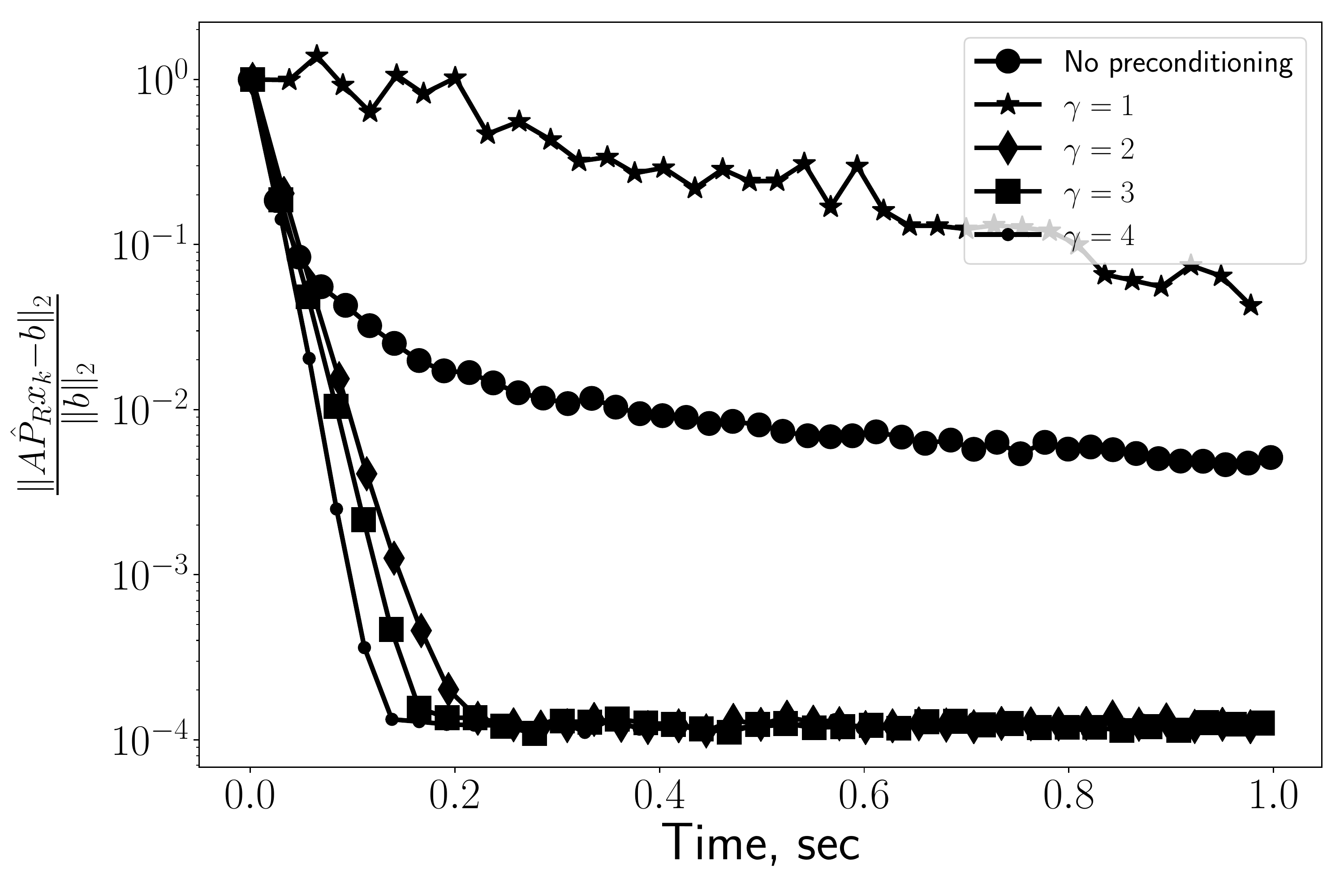}}
    \caption{Comparison of convergence speed of the considered methods for consistent linear system with noisy right-hand side for different level of noise.}
    \label{fig::random_inconsistent_time}
\end{figure}

In this setting, we generate the right-hand side $b = Ax + \varepsilon$, where $\varepsilon \in \mathcal{N}(0, \sigma^2I)$.
We test the dependence of the relative error on the required time for different values of $\gamma$ and $\sigma$.
Figure~\ref{fig::random_inconsistent_time} shows the same result with respect to the considered values of $\gamma$ as for the consistent setting, see Figure~\ref{fig::random_consistent_time}.
In particular, $\gamma = 1$ is not large enough to give a good approximation of $P^*_R$, but $\gamma = 2$ or larger is already large enough to get faster convergence.
Also, this plot shows the dependence of the performance of the considered methods for different noise levels.
The higher the level of noise, the less gain from the preconditioned Kaczmarz method we get.
That is because preconditioning gives a very small relative error, but in the high noise case, we can not get very accurate $x_k$ due to the noisy right-hand side.

% \begin{figure}[!h]
%     \centering
%     \begin{subfigure}[t]{0.35\textwidth}
%     \centering
%     \includegraphics[scale=0.25]{{inconsistent_sigma_0.1_time}.pdf}
%     \caption{$\sigma = 10^{-1}$}
%     \end{subfigure}
%     ~
%     \begin{subfigure}[t]{0.35\textwidth}
%     \centering
%     \includegraphics[scale=0.25]{{inconsistent_sigma_0.01_time}.pdf}
%     \caption{$\sigma = 10^{-2}$}
%     \end{subfigure}
%     \\
%     \begin{subfigure}[t]{0.47\textwidth}
%     \includegraphics[scale=0.25]{{inconsistent_sigma_0.001_time}.pdf}
%     \caption{$\sigma = 10^{-3}$}
%     \end{subfigure}
%     ~
%     \begin{subfigure}[t]{0.47\textwidth}
%     \includegraphics[scale=0.25]{{inconsistent_sigma_0.0001_time}.pdf}
%     \caption{$\sigma = 10^{-4}$}
%     \end{subfigure}
%     \caption{Comparison of convergence speed of the considered methods for consistent linear system with noise for different level of noise.}
%     \label{fig::random_inconsistent_time}
% \end{figure}

Thus, we experimentally show the effectiveness of the preconditioning Kaczmarz method in the cases of randomly generated consistent linear system and randomly generated linear system with a small level of noise in the right-hand side.
In the next section, we use matrices coming from different types of tomography to show the performance of the proposed method in the real application. 

\subsection{Tomography data}
\begin{figure}[!h]
    \centering
    \includegraphics[scale=0.5]{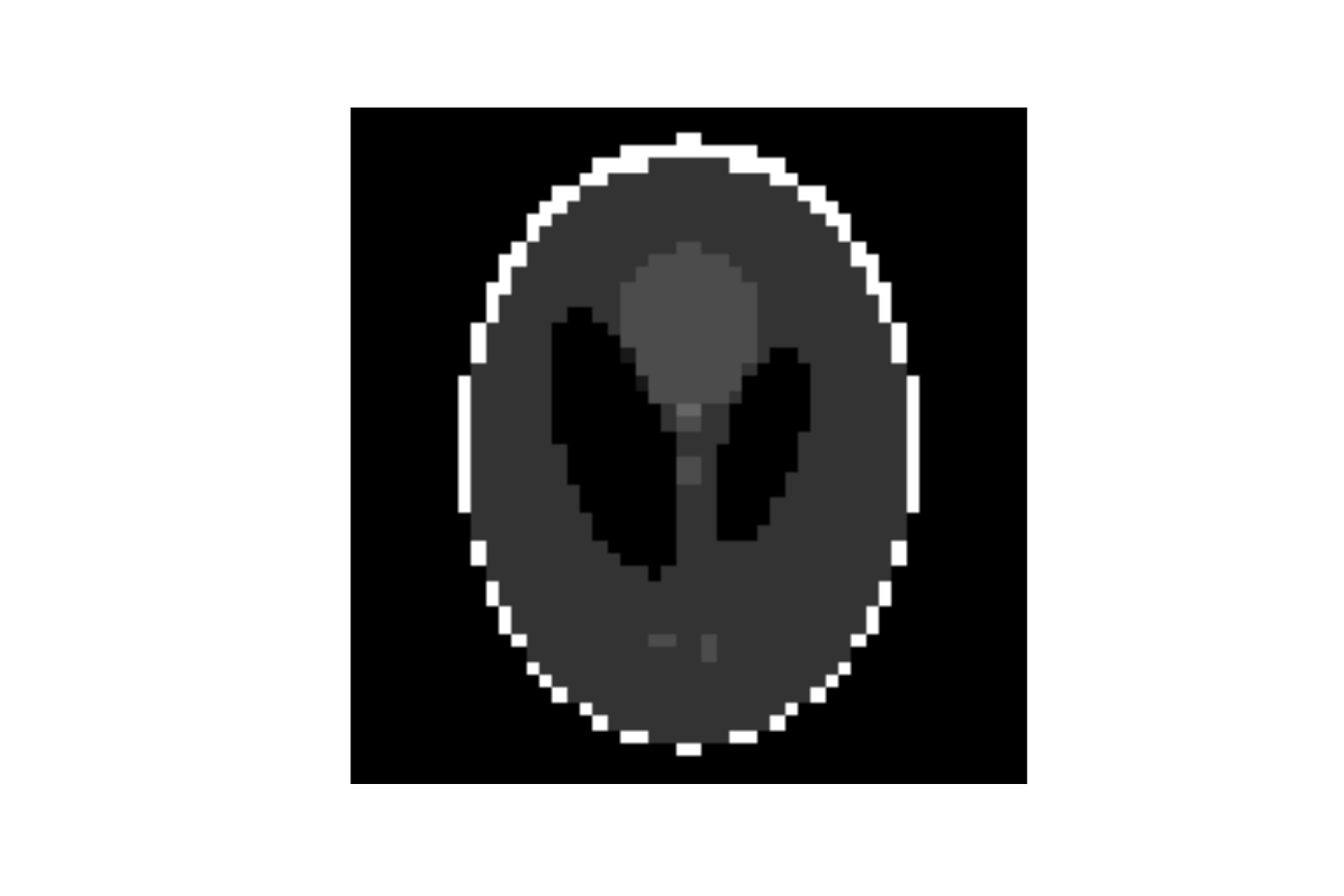}
    \caption{Original test image}
    \label{fig::shepplogan}
\end{figure}

To generate the tomography data we use the package AIR Tool II~\cite{hansen2018air}.
In particular, to generate matrices $A$ we use functions \texttt{paralleltomo}, \texttt{fancurvedtomo} and \texttt{fanlineartomo} with default parameters and number of discretization intervals $q = 50$.
The number of intervals $q$ corresponds to a $q \times q$ square image that we reconstruct through solving an overdetermined system.
With default parameters and $q = 50$ these functions give matrices $12780 \times 2500$.
Denote by $A_p, A_c$ and $A_l$ matrices given by functions \texttt{paralleltomo}, \texttt{fancurvedtomo} and \texttt{fanlineartomo} respectively.
As the ground truth solution we use a standard test image shown in Figure~\ref{fig::shepplogan} and denote it by $X^* \in \mathbb{R}^{q \times q}$.
This image is represented as $q \times q$ matrix and is reconstructed by the considered methods in the reshaped form $q^2 \times 1$.

\begin{figure}[!h]
    \centering
    \subfloat[Matrix $A_p$]{\label{fig::parallel_conv_shepplogan} \includegraphics[scale=0.21]{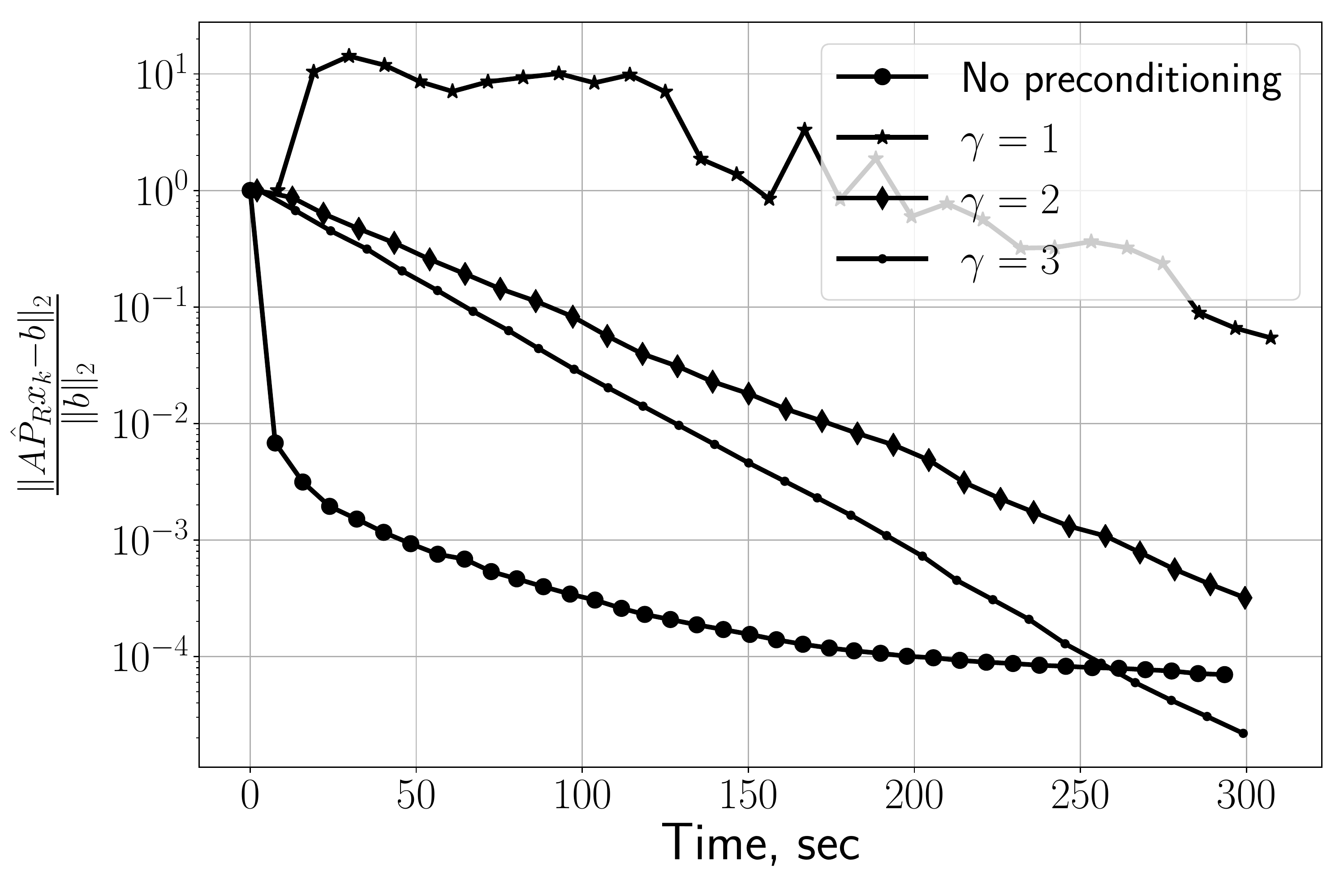}}
    \subfloat[Matrix $A_c$]{\label{fig::curved_conv_shepplogan}\includegraphics[scale=0.21]{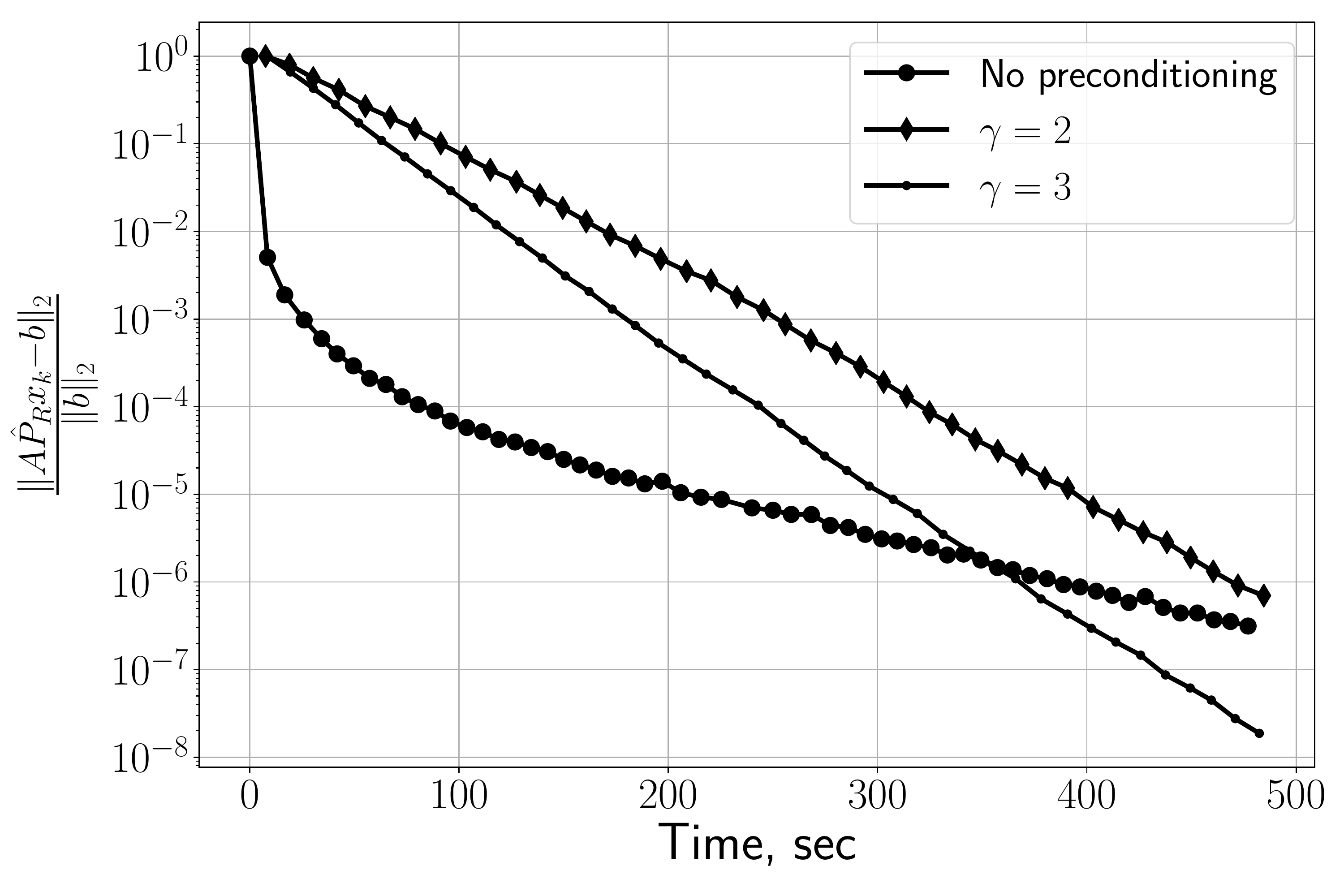}}
    \\
    \subfloat[Matrix $A_l$]{\label{fig::linear_conv_shepplogan} \includegraphics[scale=0.21]{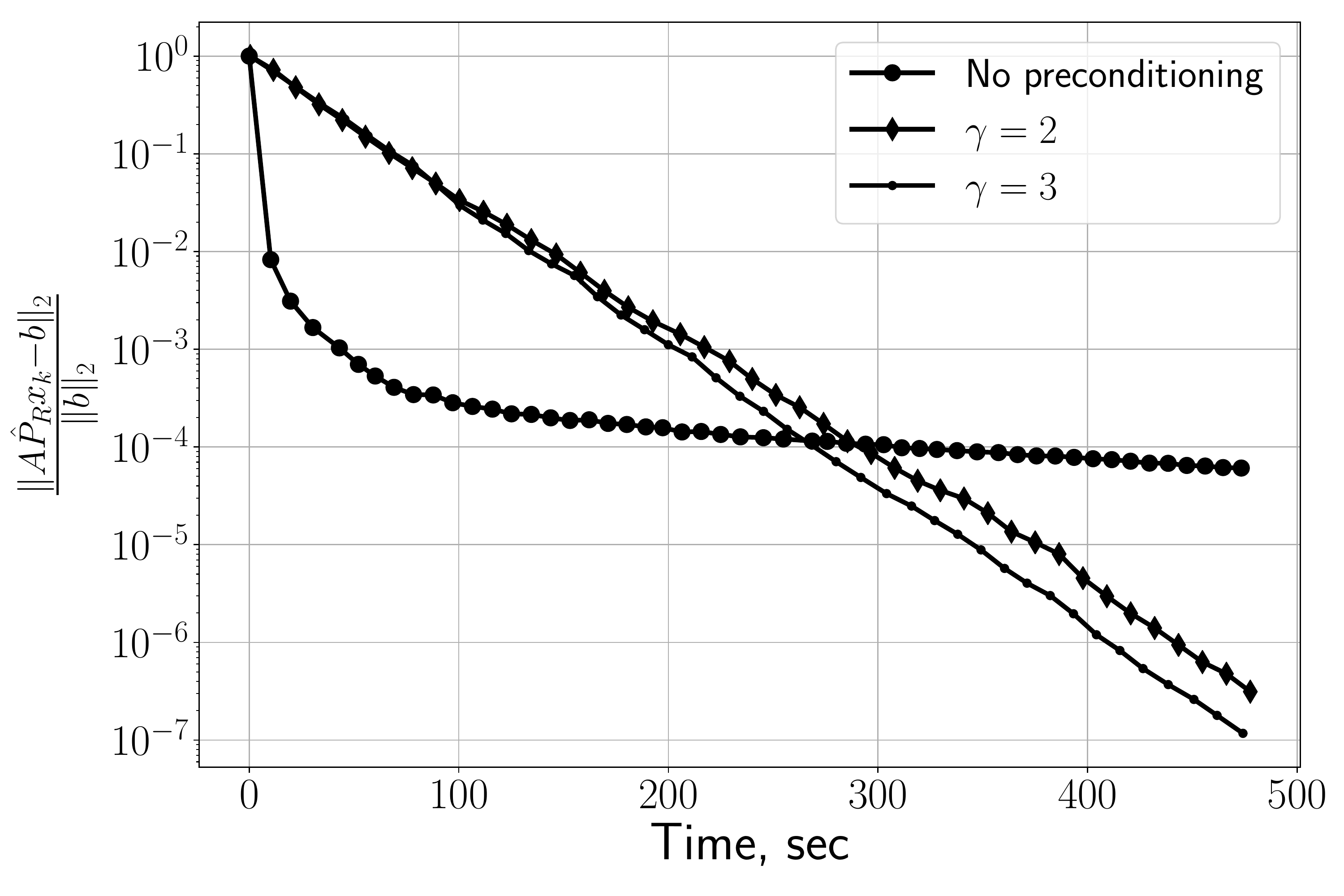}}
    \caption{Convergence comparison of the considered methods for the test image, see Figure~\ref{fig::shepplogan}}
    \label{fig::conv_shepplogan}
\end{figure}

Figure~\ref{fig::parallel_conv_shepplogan} shows the convergence of the same methods that were studied in Section~\ref{sec::random_data}.
Here we again observe that $\gamma = 1$ is not enough for faster convergence.
Moreover, for the matrix $A_p$, $\gamma = 1$ gives a singular matrix $\hat{R}$ and we use its pseudoinverse to compute $\hat{P}_R$.
For larger values of $\gamma$ we get nonsingular matrices $\hat{R}$ and observe asymptotically faster convergence compared with no preconditioning case.
However, we see two phases in the convergence of preconditioned methods: the first phase gives less accurate approximation than the Kaczmarz method without preconditioning, but at some moment the Kaczmarz method gives higher relative error than the preconditioned Kaczmarz methods.
For example, Figure~\ref{fig::parallel_conv_shepplogan} shows that after 300 seconds the preconditioned Kaczmarz method, where $\gamma = 3$, gives lower relative error than the Kaczmarz method without preconditioning.
We observe the same convergence in the case of matrices $A_c$ and $A_l$, see Figure~\ref{fig::curved_conv_shepplogan},~\ref{fig::linear_conv_shepplogan}.
To make the convergence plot more clear, we skip the case $\gamma = 1$ since it gives non-singular but very ill-conditioned matrix $\hat{R}$.

% \begin{figure}[!h]
%     \centering
%     \begin{subfigure}[t]{0.45\textwidth}
%     \centering
%     \includegraphics[scale=0.25]{parallel_time_shepplogan.pdf}
%     \caption{Matrix $A_p$}
%     \label{fig::parallel_conv_shepplogan}
%     \end{subfigure}
%     ~
%     \begin{subfigure}[t]{0.45\textwidth}
%     \centering
%     \includegraphics[scale=0.25]{curved_time_shepplogan.pdf}
%     \caption{Matrix $A_c$}
%     \label{fig::curved_conv_shepplogan}
%     \end{subfigure}
%     \\
%     \begin{subfigure}[t]{0.45\textwidth}
%     \centering
%     \includegraphics[scale=0.25]{linear_time_shepplogan.pdf}
%     \caption{Matrix $A_l$}
%     \label{fig::linear_conv_shepplogan}
%     \end{subfigure}
%     \caption{Convergence comparison of the considered methods for the test image, see Figure~\ref{fig::shepplogan}}
%     \label{fig::conv_shepplogan}
% \end{figure}

\begin{figure}[!h]
    \subfloat[Matrix $A_p$]{\includegraphics[scale=0.15]{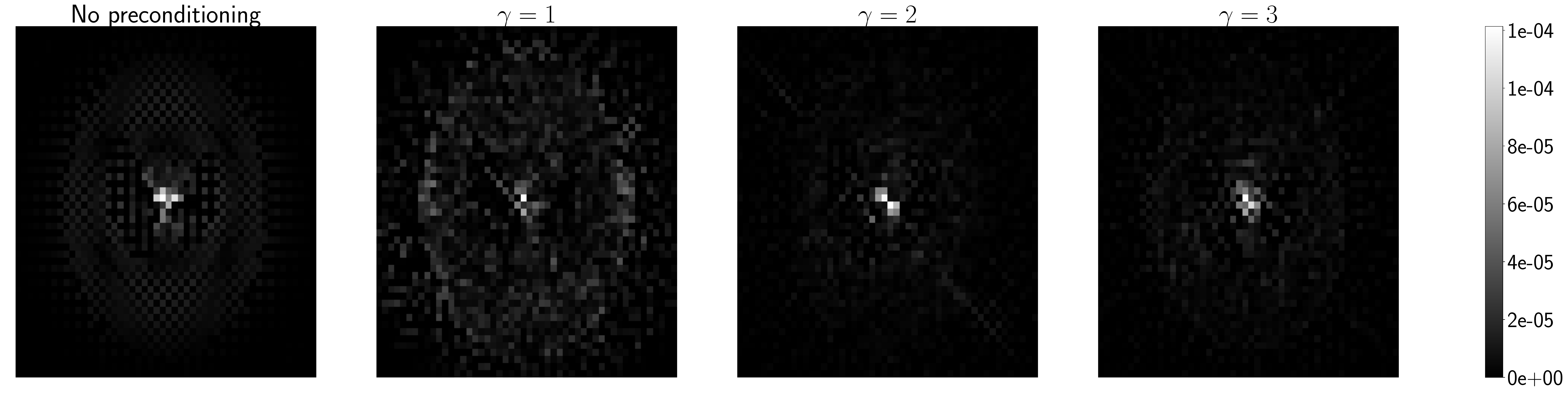}}
    \\
    \subfloat[Matrix $A_c$]{\includegraphics[scale=0.19]{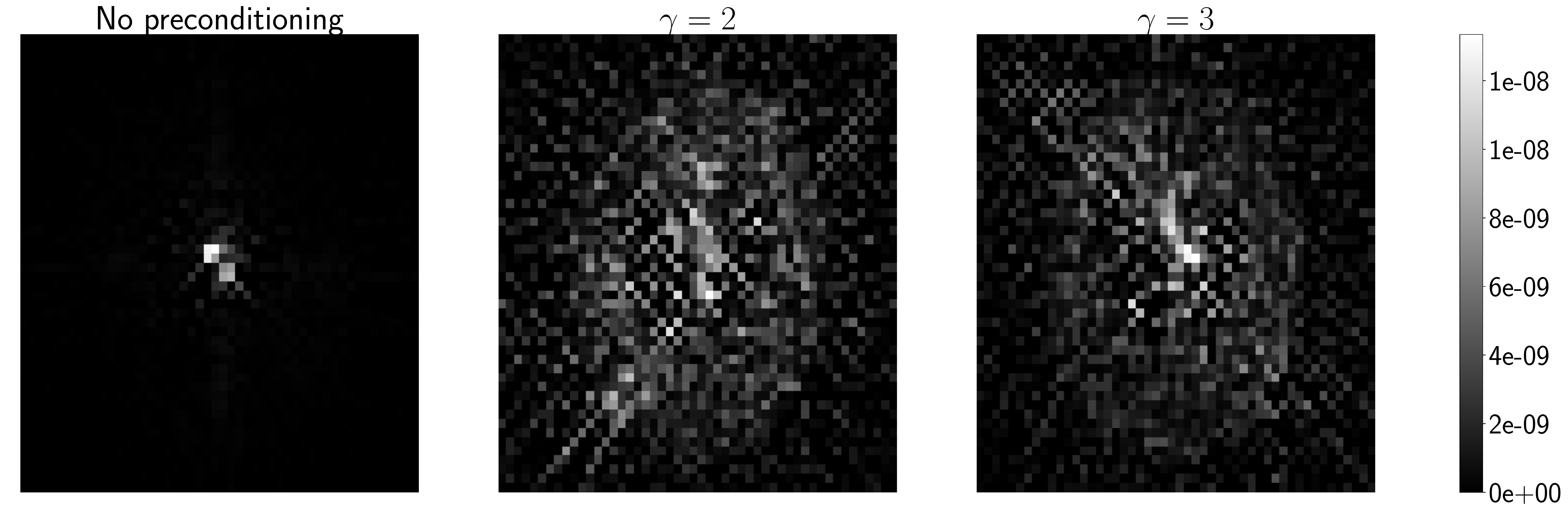}}
    \\
    \subfloat[Matrix $A_l$]{\includegraphics[scale=0.19]{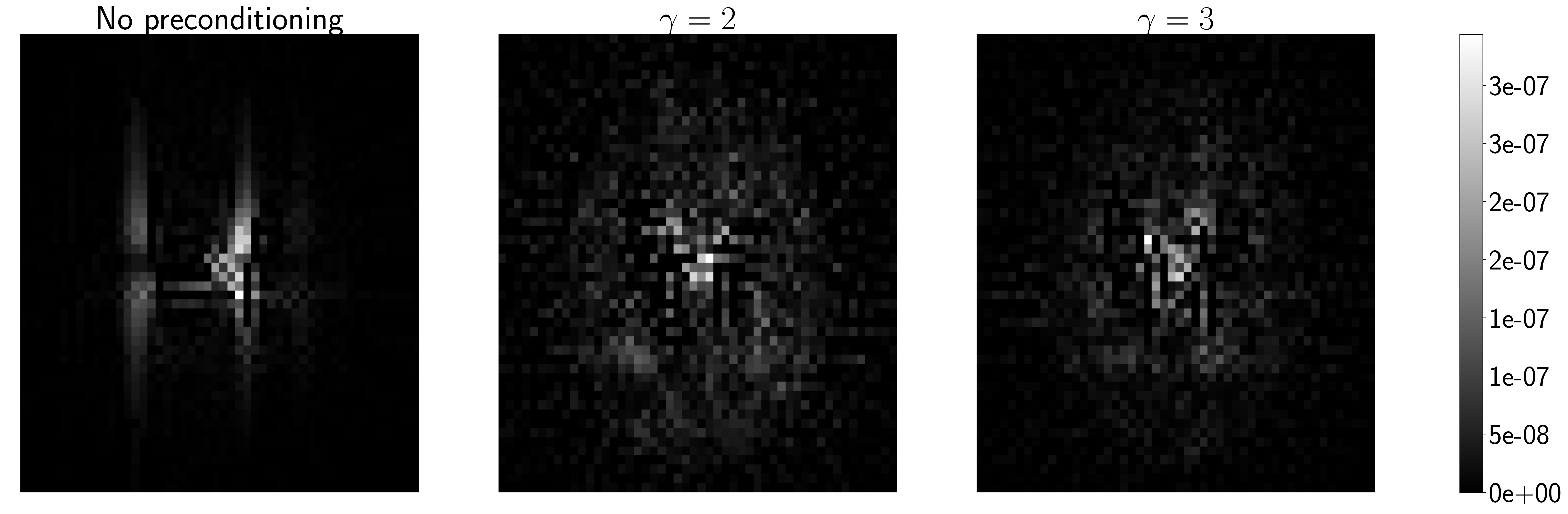}}
    \caption{Relative error for different tomography types computed as $\frac{|X^* - \hat{X}|}{\|X^*\|_F}$, where nominator is elementwise absolute value of difference between ground truth solution $X^*$ and approximation $\hat{X}$ given by considered methods}
    \label{fig::error_mask_shepplogan}
\end{figure}

In addition to the convergence plots in Figure~\ref{fig::conv_shepplogan}, we provide the normalized elementwise error between ground truth image and images reconstructed with the considered methods, see Figure~\ref{fig::error_mask_shepplogan}.
These plots show regions of the reconstructed images with high and low reconstruction quality, which is measured as $\frac{|X^* - \hat{X}|}{\|X^*\|_F}$, where $X^*$ is ground truth solution and $\hat{X}$ is approximation given by the corresponding method. 
The nominator is the elementwise absolute value of the difference $X^*$~and~$\hat{X}$.  

% \begin{figure}[!h]
%     \begin{subfigure}[t]{\textwidth}
%     \centering
%     \includegraphics[scale=0.2]{error_mask_parallel_shepplogan.pdf}
%     \caption{Tomography generated by matrix $A_p$}
%     \end{subfigure}
%     \\
%     \begin{subfigure}[t]{\textwidth}
%     \centering
%     \includegraphics[scale=0.2]{error_mask_curved_shepplogan.pdf}
%     \caption{Tomography generated by matrix $A_c$}
%     \end{subfigure}
%     \\
%     \begin{subfigure}[t]{\textwidth}
%     \centering
%     \includegraphics[scale=0.2]{error_mask_linear_shepplogan.pdf}
%     \caption{Tomography generated by matrix $A_l$}
%     \end{subfigure}
%     \caption{Relative error for different tomography types computed as $\frac{|X^* - \hat{X}|}{\|X^*\|_F}$, where nominator is elementwise absolute value of difference between ground truth solution $X^*$ and approximation $\hat{X}$ given by considered methods}
%     \label{fig::error_mask_shepplogan}
% \end{figure}

Thus, the proposed preconditioned Kaczmarz method can be used as a fine-tuning technique that will improve the accuracy given by the Kaczmarz method without preconditioning.
Details of this approach are in the next section.

\subsubsection{Fine-tuning with the preconditioned Kaczmarz method}

\begin{figure}[!h]
    \centering
    \subfloat[Matrix $A_p$]{\label{fig::parallel_conv_fine_tuning_shepplogan} \includegraphics[scale=0.2]{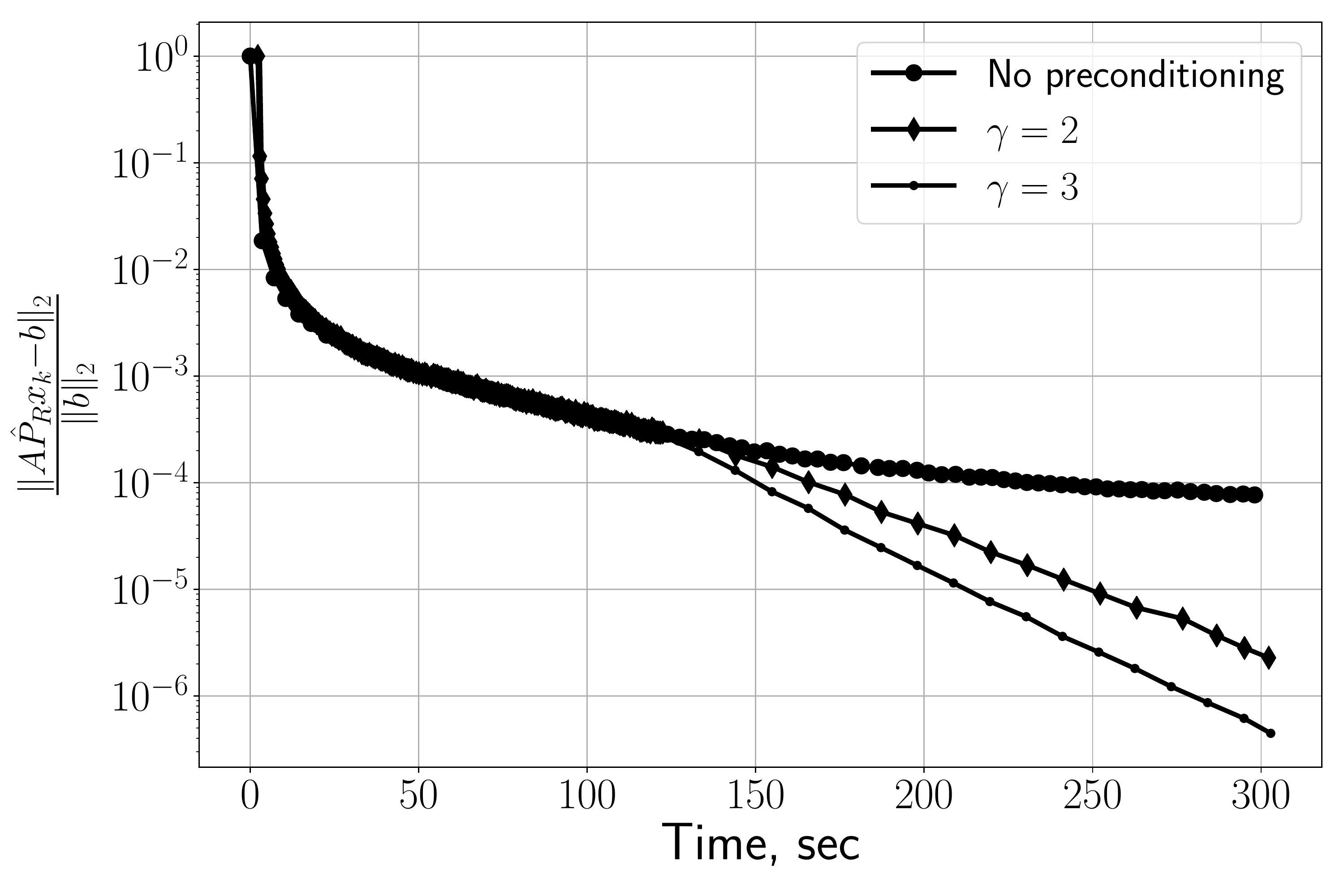}}
    \subfloat[Matrix $A_c$]{\label{fig::curved_conv_fine_tuning_shepplogan} \includegraphics[scale=0.2]{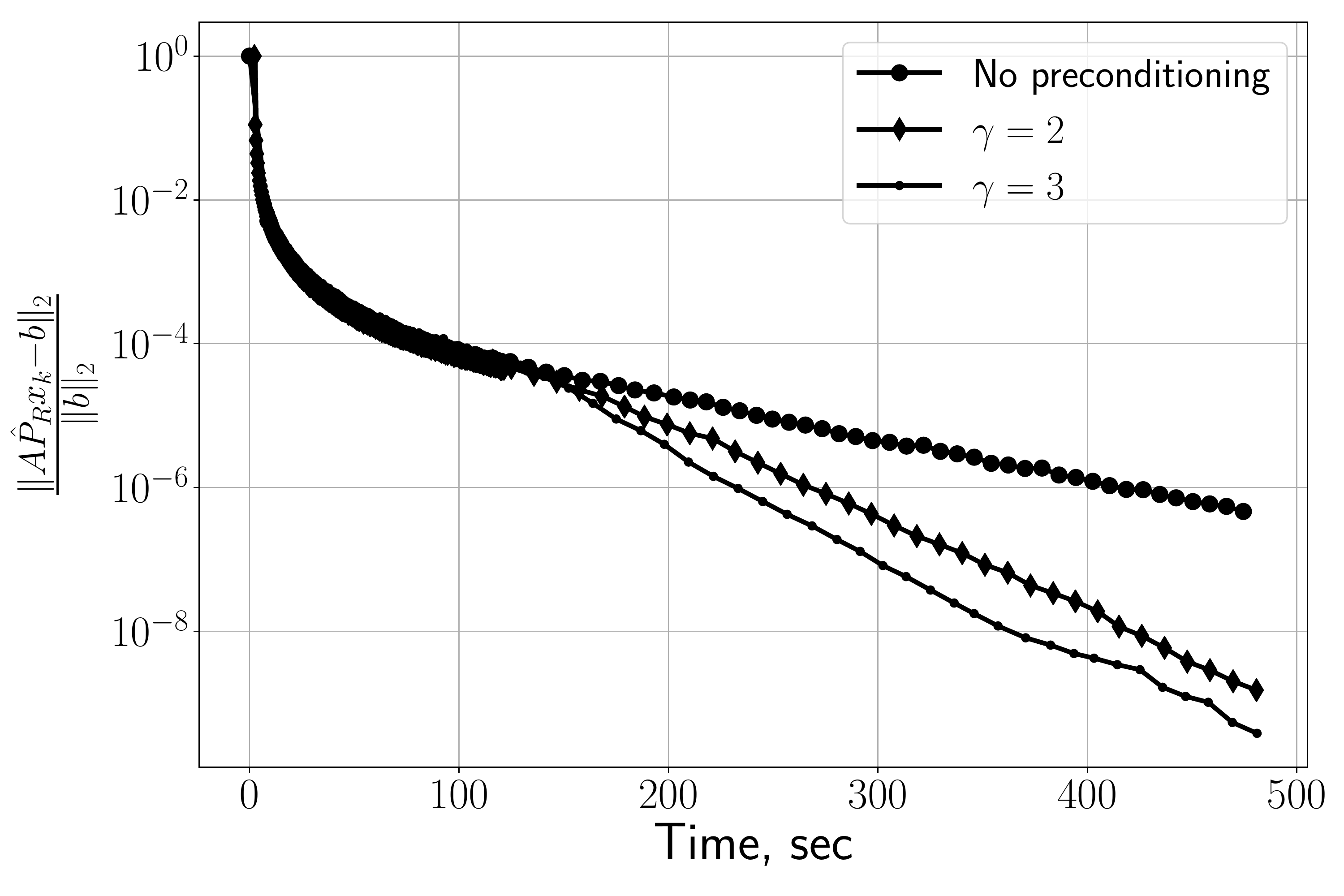}}
    \\
    \subfloat[Matrix $A_l$]{\label{fig::linear_conv_fine_tuning_shepplogan} \includegraphics[scale=0.2]{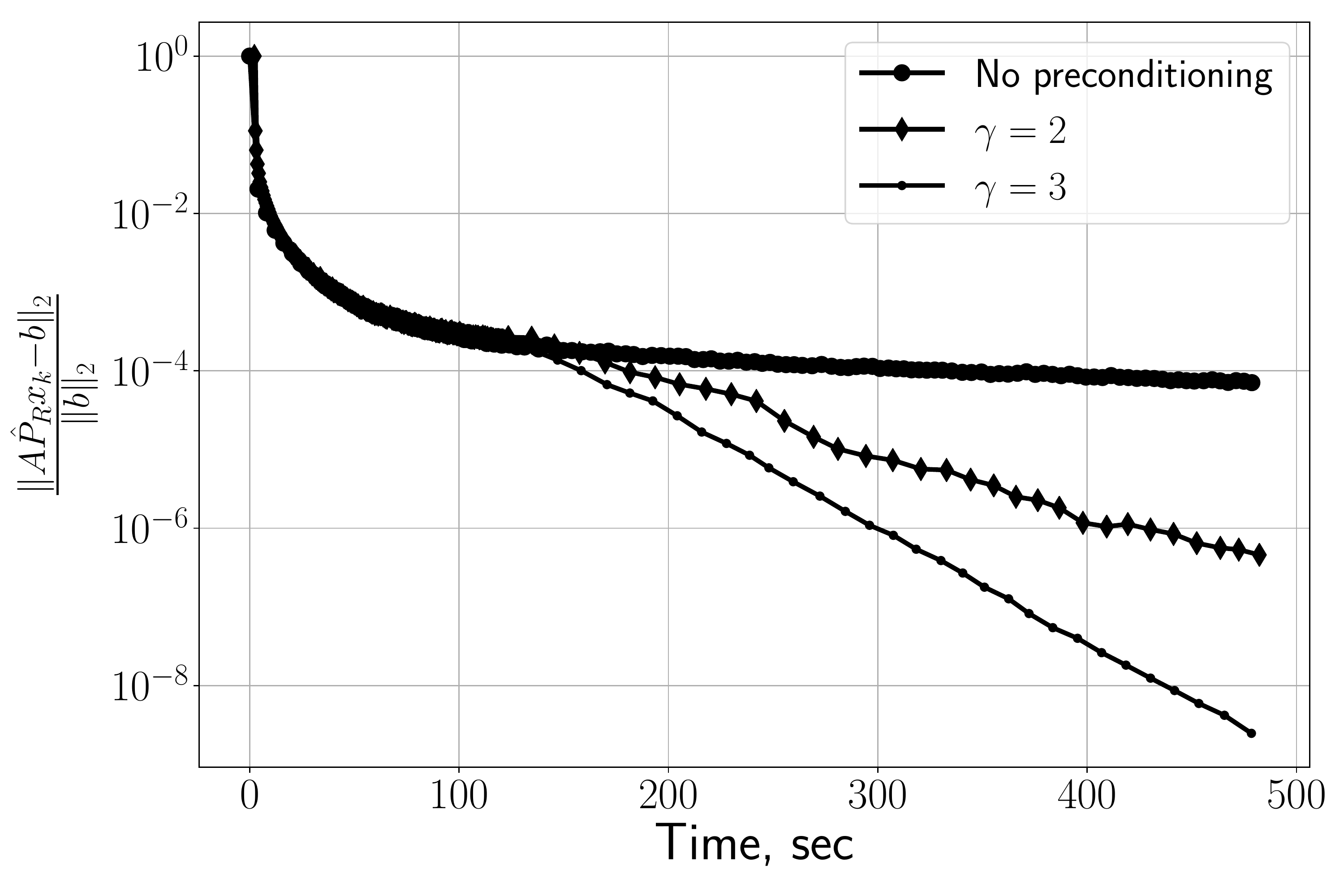}}
    \caption{Convergence comparison of the considered methods with fine-tuning strategy ($\tau = 120$) for the test image, see Figure~\ref{fig::shepplogan}}
    \label{fig::conv_fine_tuning_shepplogan}
\end{figure}

As we show in the previous section the preconditioned Kaczmarz asymptotically gives a smaller relative error. 
However, the Kaczmarz method without preconditioning gives a smaller relative error during the first iterations. 
Therefore, we propose to start solving a given system without using preconditioner and then switch to the preconditioned modification of the Kaczmarz method.
This strategy gives a lower error rate for the same runtime, see Figure~\ref{fig::conv_fine_tuning_shepplogan}, where we denote by $\tau$ the time to switch to the  preconditioned Kaczmarz method.
We use $\tau = 120$ seconds in all experiments.
In the case of tomography given by $A_p$ we do not consider fine-tuning for $\gamma = 1$, since it gives poor convergence that we see in Figure~\ref{fig::parallel_conv_shepplogan}.

% \begin{figure}[!h]
%     \centering
%     \begin{subfigure}[t]{0.45\textwidth}
%     \centering
%     \includegraphics[scale=0.25]{parallel_time_fine_tuning120_shepplogan.pdf}
%     \caption{Matrix $A_p$}
%     \label{fig::parallel_conv_fine_tuning_shepplogan}
%     \end{subfigure}
%     ~
%     \begin{subfigure}[t]{0.45\textwidth}
%     \centering
%     \includegraphics[scale=0.25]{curved_time_fine_tuning120_shepplogan.pdf}
%     \caption{Matrix $A_c$}
%     \label{fig::curved_conv_fine_tuning_shepplogan}
%     \end{subfigure}
%     \\
%     \begin{subfigure}[t]{0.45\textwidth}
%     \centering
%     \includegraphics[scale=0.25]{linear_time_fine_tuning120_shepplogan.pdf}
%     \caption{Matrix $A_l$}
%     \label{fig::linear_conv_fine_tuning_shepplogan}
%     \end{subfigure}
%     \caption{Convergence comparison of the considered methods with fine-tuning strategy ($\tau = 120$) for the test image, see Figure~\ref{fig::shepplogan}}
%     \label{fig::conv_fine_tuning_shepplogan}
% \end{figure}

\subsection{Kernel regression problem}
One more source of the consistent overdetermined linear systems is the kernel regression problem with  explicit feature maps, e.g. \emph{random Fourier features}~\cite{avron2017random,munkhoeva2018quadrature}.
These maps reduce the non-linear regression problem to the linear regression problem with the modified matrix.
We use this approach in the problem of non-linear function approximation such that the parametric function model is not known.
Also, it is known that introducing Gaussian kernel reduces linear regression problem to non-linear and can give a very accurate approximation of smooth functions~\cite{micchelli2006universal}.

Consider the following function $f: \mathbb{R}^2 \to \mathbb{R}$ that is used as the test function in the package George~\cite{ambikasaran2016fast}:
\begin{equation}
f(x) = a^{\top}x + c + \alpha \exp\left( -\frac{(\|x\|_2 - \mu)^2}{\sigma_f}\right),
\label{eq::test_func_kernel}
\end{equation}
and use the following parameters: $a = (1, 1), c = 1, \alpha=0.1, \mu=0, \sigma_f=1$.
To compose an overdetermined linear system we sample $m = 10000$ points uniform in $D = [0, 1] \times [0, 1]$ and compute the function values in these points.
So, we have a matrix~$Z \in \mathbb{R}^{m \times 2}$, where every row is a sampled point from the set~$D$ and the right-hand side~$b$ is the vector of the function values in the corresponding points.
Since the test function~\eqref{eq::test_func_kernel} is non-linear, the composed linear system~$Z\tilde{x} = b$ is not consistent, but the implicit kernel trick with random Fourier features helps in constructing a \emph{consistent} linear system $Ax = b$.
This explicit feature map gives an approximation of introducing the Gaussian kernel and the accuracy of the approximation exponentially grows with the dimension~$d$~\cite{rahimi2008random}. 
To construct the matrix $A$, we generate the matrix $M \in \mathbb{R}^{2 \times d}$ with entries from the normal distribution $\mathcal{N}(0, \sigma^2)$, where $d \geq 2$ and $\sigma$ are hyper-parameters that have to be tuned.
After that, we compute $B = ZM \in \mathbb{R}^{m \times d}$, $B_c = \cos(B)$ and $B_s = \sin(B)$, where sine and cosine are computed elementwise.
Now we are ready to construct the matrix $A \in \mathbb{R}^{m \times 2d}$ that gives a consistent linear system. 
We use the columns of the matrices $B_c$ and $B_s$ in the following way:
\begin{enumerate}
    \item[1)] the $i$-th column of the matrix $B_c$ is the $(2i-1)$-th column of the matrix $A$ for $i=1,\ldots, d$
    \item[2)] the $i$-th column of the matrix $B_s$ is the $2i$-th column of the matrix $A$ for $i=1,\ldots, d$.
\end{enumerate}
Thus, for sufficiently large $d$ and small $\sigma$ we have a consistent linear system $Ax = b$~\cite{avron2017random}.

\begin{figure}[!h]
    \centering
    \subfloat[$d = 3$]{\includegraphics[scale=0.2]{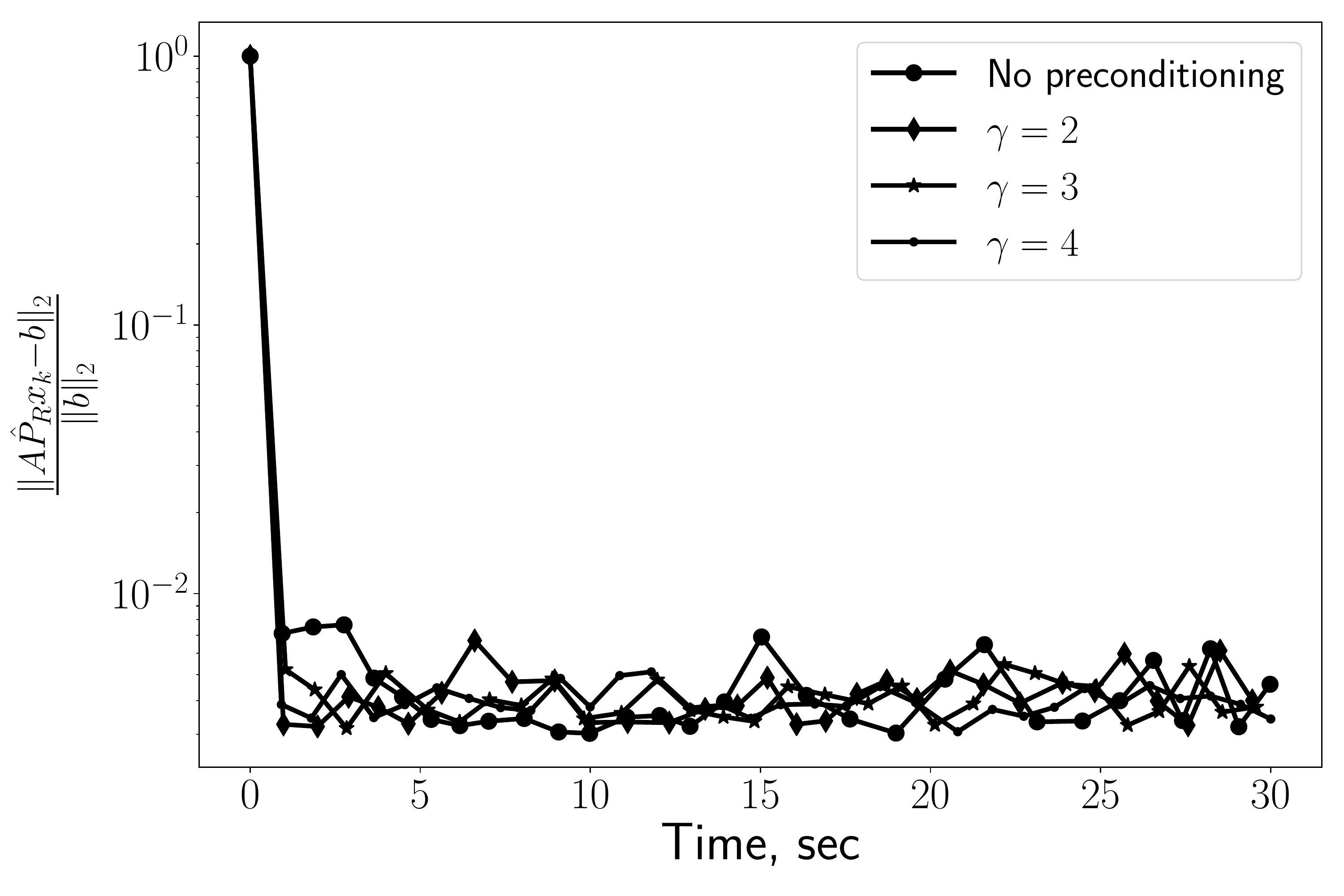}}
    \subfloat[$d = 5$]{\label{fig::rff_d5}\includegraphics[scale=0.2]{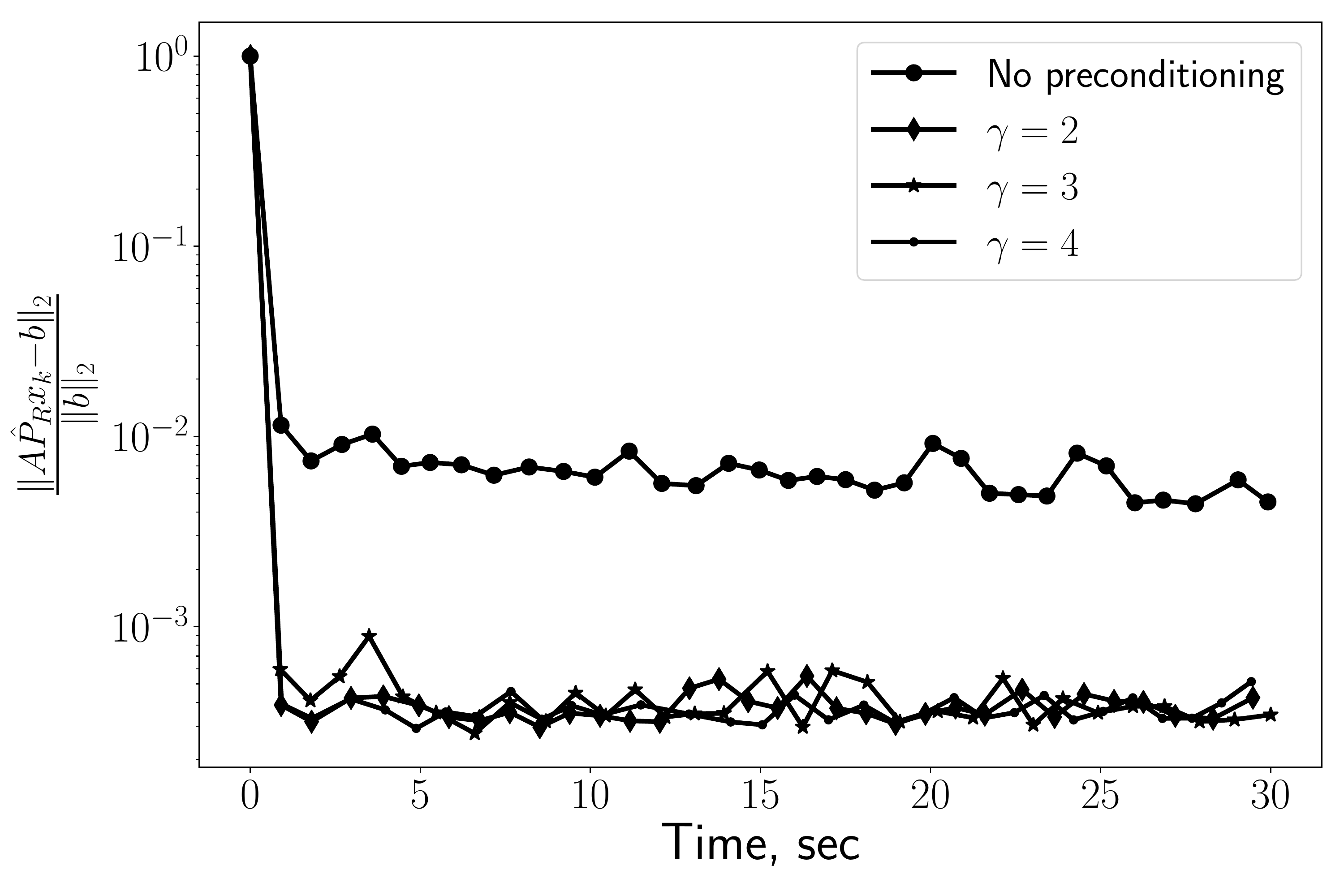}}
    \\
    \subfloat[$d = 10$]{\includegraphics[scale=0.2]{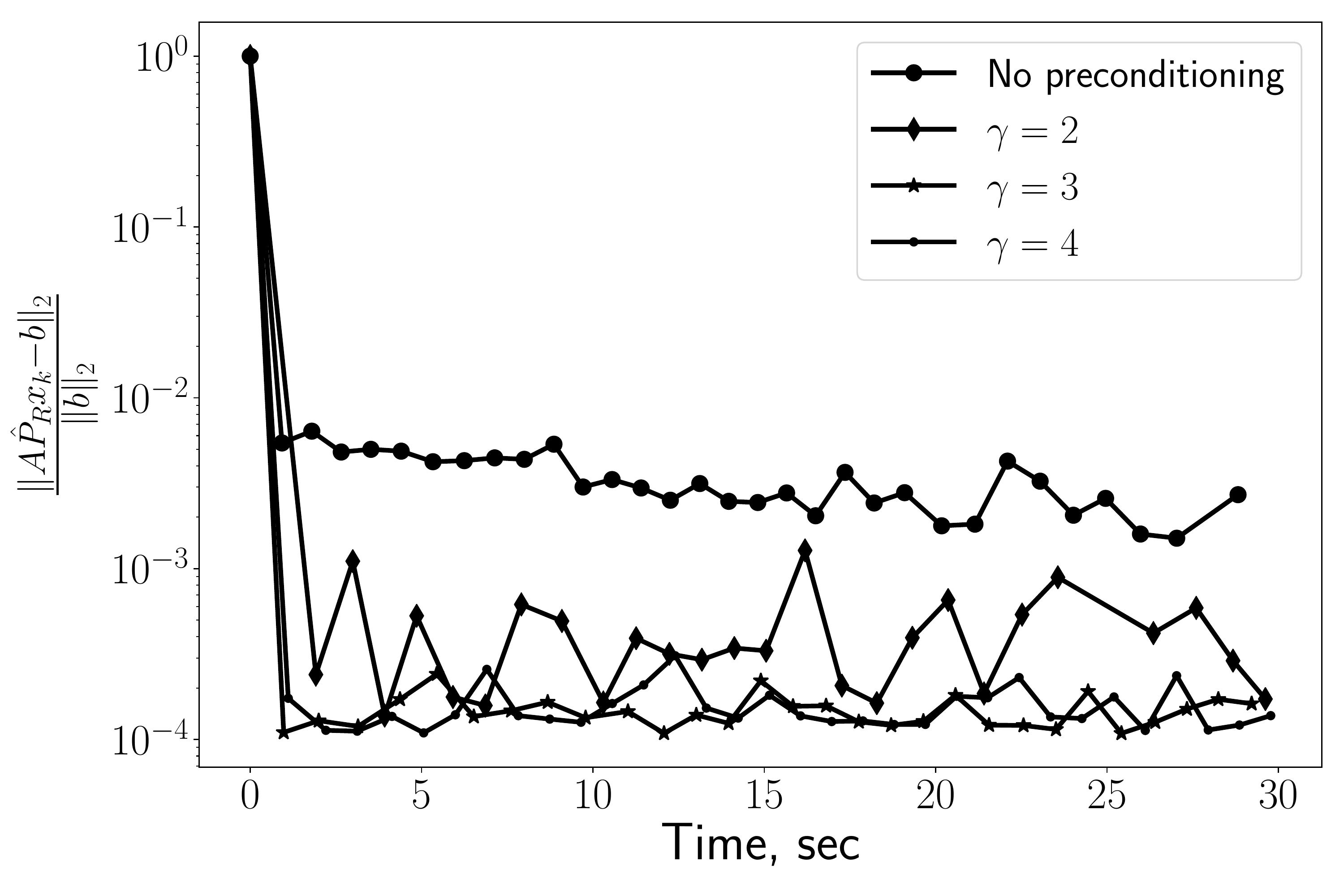}}
    \subfloat[$d = 20$]{\includegraphics[scale=0.2]{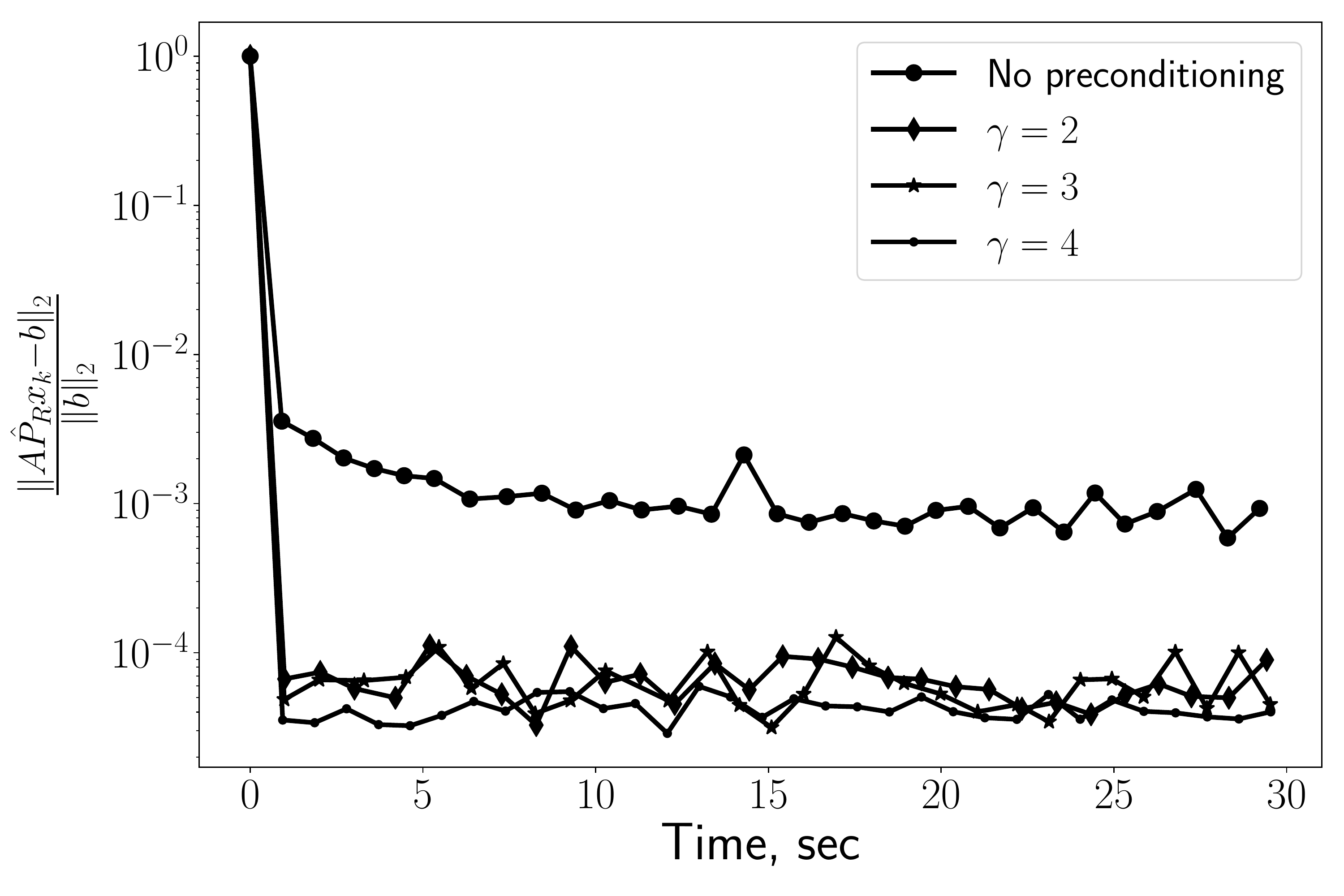}}
    \caption{Convergence comparison of the considered methods for the kernel regression problem via random Fourier features for function~\eqref{eq::test_func_kernel} approximation, $\sigma=1$}
    \label{fig::rff}
\end{figure}

Figure~\ref{fig::rff} shows the convergence of the considered methods for different values of $d$.
Experiments show that $\sigma = 1$ gives a sufficiently accurate solution, therefore we fix it and show the dependence only on the dimension $d$. 
Figure~\ref{fig::rff_d5} shows that already for $d = 5$ the preconditioned Kaczmarz method converges much faster for all considered values of $\gamma$.

\section{Conclusion}
This paper presents the preconditioned Kaczmarz method by sketching technique and demonstrates its performance in different applications.
However, the proposed approach requires a number of sketched rows that is an unknown hyperparameter, therefore theoretical investigation on the minimal, sufficient number of the sketched rows is the interesting further study.
Another ingredient of the Kaczmarz method is the distribution for a sampling of rows. 
To sample rows in the right way, the methods for sampling from unnormalized distributions can be used.
Also, the preconditioning approach can be generalized to non-linear systems and inspire modifications of the corresponding non-linear stochastic least-squares solvers.
One more interesting extension of this work is to consider state-of-the-art adaptive stochastic gradient methods~\cite{duchi2011adaptive,kingma2014adam} from the preconditioning perspective like in~\cite{staib2019escaping}. 
This point of view can lead to more efficient and theoretically motivated optimization methods.

% In experiments we use synthetic random consistent and with noisy right-hand side linear systems. 
% Also, to test the considered methods we use tomography data and non-linear function approximation problem.
% This problem can be reduced to consistent overdetermined linear system by explicit feature map approach. 
% We use random Fourier feature as such a map that approximates explicit introducing Gaussian kernel function.
% Experimental results show that the proposed method is faster than the Kaczmarz method without preconditioning.
% Also we study linear systems with a noisy right-hand side and observe that the higher level of noise, the less gain is obtained.
% We find that the most effective way to use the preconditioned Kaczmarz method in the case of large number of columns is fine-tuning.
% The Kaczmarz method gives less error rate than preconditioned one for the first iterations.
% Therefore, after these first iterations we switch to the preconditioned Kaczmarz method and get high accurate approximation.

\bibliographystyle{siamplain}
\bibliography{lib}

\end{document}